\documentclass{elsart}

\usepackage{graphicx}

\usepackage{amssymb}

\begin{document}

\begin{frontmatter}



\title{The obstructions for toroidal graphs with no $K_{3,3}$'s}


\author[Montreal]{Andrei Gagarin\corauthref{cor}}
\corauth[cor]{Corresponding author.}

\address[Montreal]{Laboratoire de Combinatoire et d'Informatique Math\'ematique (LaCIM), PK-4211,
Universit\'e du Qu\'ebec \`a Montr\'eal, C.P. 8888, Succ. Centre-Ville, Montr\'eal, Qu\'ebec, H3C 3P8, Canada}
\ead{gagarin@lacim.uqam.ca}
\ead[url]{http://www.lacim.uqam.ca/$\sim$gagarin/}

\author[Victoria]{Wendy Myrvold\thanksref{NSERC}}, and 
\ead{wendym@cs.uvic.ca}
\ead[url]{http://csr.uvic.ca/$\sim$wendym/}
\author[Victoria]{John Chambers}
\address[Victoria]{Department of Computer Science, University of Victoria, P.O. Box 3055, Stn. CSC, 
Victoria, British Columbia, 
V8W 3P6, Canada}
\thanks[NSERC]{This work was supported by an operating grant from the Natural Sciences and Engineering Research Council of Canada.}

\begin{abstract}
Forbidden minors and subdivisions for toroidal graphs are numerous. We consider the toroidal graphs with no $K_{3,3}$-subdivisions that coincide with the toroidal graphs with no $K_{3,3}$-minors. These graphs admit a unique decomposition into planar components and have short lists of obstructions. We provide the complete lists of four forbidden minors and eleven forbidden subdivisions for the toroidal graphs with no $K_{3,3}$'s and prove that the lists are sufficient.
\end{abstract}

\end{frontmatter}

\section{Introduction}
\noindent We use basic graph-theoretic terminology from Bondy and Murty \cite{Bondy} and Diestel \cite{Diestel}. An {\it undirected graph} $G$ consists of a set $V$ of {\it vertices} and a set $E$ of {\it edges} where each edge is associated with an unordered pair $uv$ of vertices in $V$. A graph $G$ is {\it embeddable} in a surface if it can be drawn on the surface with no crossing edges. A {\it topological obstruction} for a surface is a graph $G$ with minimum vertex degree three which does not embed on the surface, but $G-e$ is embeddable for any edge $e$. A {\it minor order obstruction} is a topological obstruction with the additional property that $G$ with edge $e$ contracted embeds on the surface for all edges $e$.

The purpose of this paper is to characterize all the minor order and topological obstructions for the torus which do not contain a subgraph homeomorphic ({\it homeomorphic} means the edges in the figure correspond to paths in the graph) to $K_{3,3}$ pictured in Figure~1. If $H$ is a graph with minimum vertex degree three, then following Diestel \cite{Diestel}, an {\it H-subdivision} denoted by $TH$ is a graph homeomorphic to $H$.

Kuratowski's theorem \cite{Kuratowski} states that a graph $G$ is non-planar
if and only if it contains a subgraph homeomorphic to $K_{3,3}$ or $K_5$ (see Figure~1). Wagner \cite{Wagner2} proved that a graph $G$ is non-planar if and only if it contains $K_{3,3}$ or $K_5$ as a minor. We concentrate on the graphs with no $K_{3,3}$-subdivisions.
Since $K_{3,3}$ is 3-regular, it is possible to see that a graph does not contain any $K_{3,3}$-subdivision if and only if it does not contain any $K_{3,3}$-minor. Therefore we refer to these graphs as to graphs {\it with no $K_{3,3}$'s}.
\begin{figure}[h]
	\centerline {\includegraphics[width=2.7in]{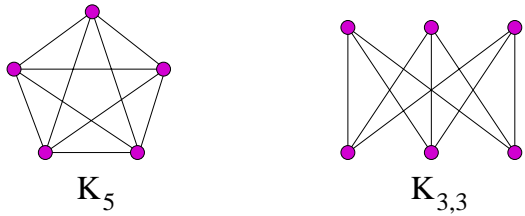}}
	\caption{Minimal non-planar graphs $K_{3,3}$ and $K_5$}
\end{figure}

Wagner \cite{Wagner} and Kelmans \cite{Kelmans} gave a recursive decomposition for all non-planar graphs with no $K_{3,3}$'s. Gagarin, Labelle and Leroux \cite{GLL2} explicitly describe the structure of toroidal graphs with no $K_{3,3}$'s in terms of 2-pole planar networks substituted into the edges of non-planar core graphs. The 2-pole planar networks arise from 2-connected planar graphs. 

For any surface of fixed genus, Robertson and Seymour theory in \cite{RS90} indicates
that the number of obstructions (topological or minor order) is finite.
A proof specifically for orientable surfaces has been independently provided by Vollmerhaus \cite{Vollmer} and Bodendiek and Wagner \cite{orient} and for nonorientable surfaces by Archdeacon and Huneke
\cite{ArchHun}. For the projective plane, there are 103 topological
obstructions and 35 minor order obstructions. Glover, Huneke and Wang first proved that these
were obstructions \cite{HG79}. Archdeacon later proved that their list was complete \cite{DA81}.

The \textit{torus} is a surface shaped like a doughnut. The precise complete lists of obstructions for the torus (both minor order and topological) are still unknown. To date, Myrvold and Chambers \cite{John} have found 239,322 topological obstructions and 16,629 minor order obstructions that include those on up to eleven vertices, the 3-regular ones on up to 24 vertices, the disconnected ones and those with a cut-vertex. Previously, only a few thousand had been determined. 

The \textit{spindle surface} is formed from the sphere by identifying two points that also can be seen as pinching a torus. The two points are usually called the \textit{north} and \textit{south poles}, and after identification it is called the pinch point. The spindle surface is a pseudosurface: no neighborhood of the pinch point is homeomorphic to an open disk.

The set of minor order and topological obstructions for the spindle surface are not known either. In general, the set of graphs that embed on a pseudosurface can be not closed under the minor order, and there are pseudosurfaces which have an infinite set of obstructions. However, graphs that embed on the spindle surface are closed under the minor order, and hence there is a finite obstruction set. The topological obstructions for cubic graphs on the spindle surface are found by Archdeacon and Bonnington \cite{AB}.

Section 2 of the present paper describes the structural and algorithmic results of Gagarin and Kocay \cite{GK} for graphs containing $K_5$-subdivisions with respect to their embeddability in the torus. In Section 3, we show the lists of four minor order and eleven topological order obstructions and prove they are sufficient for toroidal graphs with no $K_{3,3}$'s. Section 4 contains suggestions for future research on obstructions for surfaces.

\section{$K_5$-Subdivisions and toroidality}
\noindent In this section, we summarize known structural results for graphs containing a $K_5$-subdivision and the corresponding algorithmic results for the torus. A $K_5$-subdivision is denoted by $TK_5$, and a $K_{3,3}$-subdivision is denoted by $TK_{3,3}$. For a graph $H$, the vertices of the $TH$ which have degree three or more are called the {\it corner vertices\/} of the $TH$.

Let $G$ be a graph and $H$ be a subgraph of $G$. A {\it bridge} of $G$ with respect to $H$
is either 
(a) an edge of $G$ which is not in $H$, but 
both its endpoints are in $H$
(plus the endpoints), or 
(b) a connected component $C$ of the graph 
$G - V(H)$ together with edges 
connecting
a vertex in $C$ to a vertex
in $H$ (and their endpoints).
Denote by $K$ the set of corner vertices of a $TK_5$. The following theorem is implicit in the work of Asano \cite{Asano}, and has also
been stated and used by Fellows and Kaschube \cite{Fellows} and Gagarin and Kocay \cite{GK}.

\begin{thm}
Let $G$ be a $2$-connected graph with a $TK_5$.
If a bridge of the set of corner vertices $K$ of the $TK_5$ contains
three or more vertices from $K$, 
then $G$ contains a $TK_{3,3}$. Otherwise,
$G$ is decomposable into a collection of bridges of $K$
such that each bridge contains exactly two corner
vertices.
\end{thm}

Given a $TK_5$ in a graph $G$ such that all of the 
bridges of $K$ contain exactly two corner vertices, 
the {\it side component} of the $TK_5$ in $G$ corresponding to the 
pair $a, b$ of corner vertices 
has its vertex set equal to the union of the
vertex sets of the bridges containing these two corners,
and the edge set is the union of the edges of
the corresponding bridges. If $F$ is a side component
corresponding to the corners $a$ and $b$, then
the {\it augmented side component} is 
$F$ if the edge $ab$ is already in $F$ and $F + ab$ otherwise.
A side component $F$ is {\it special} if the edge $ab$ is not in $F$, $F$ is planar, but the augmented side component $F+ab$ is not planar. 

Side components of a subdivision of an $M$-graph shown in Figure~2 are defined by analogy with the side components of a $K_5$-subdivision by considering pairs of adjacent vertices of the $M$-graph. The edge $xy$ of an $M$-graph with both end-points of degree seven is called the {\it central edge} of the $M$-graph (see Figure~2). 
\smallbreak
\begin{figure}[h]
	\centerline {\includegraphics[width=1.5in]{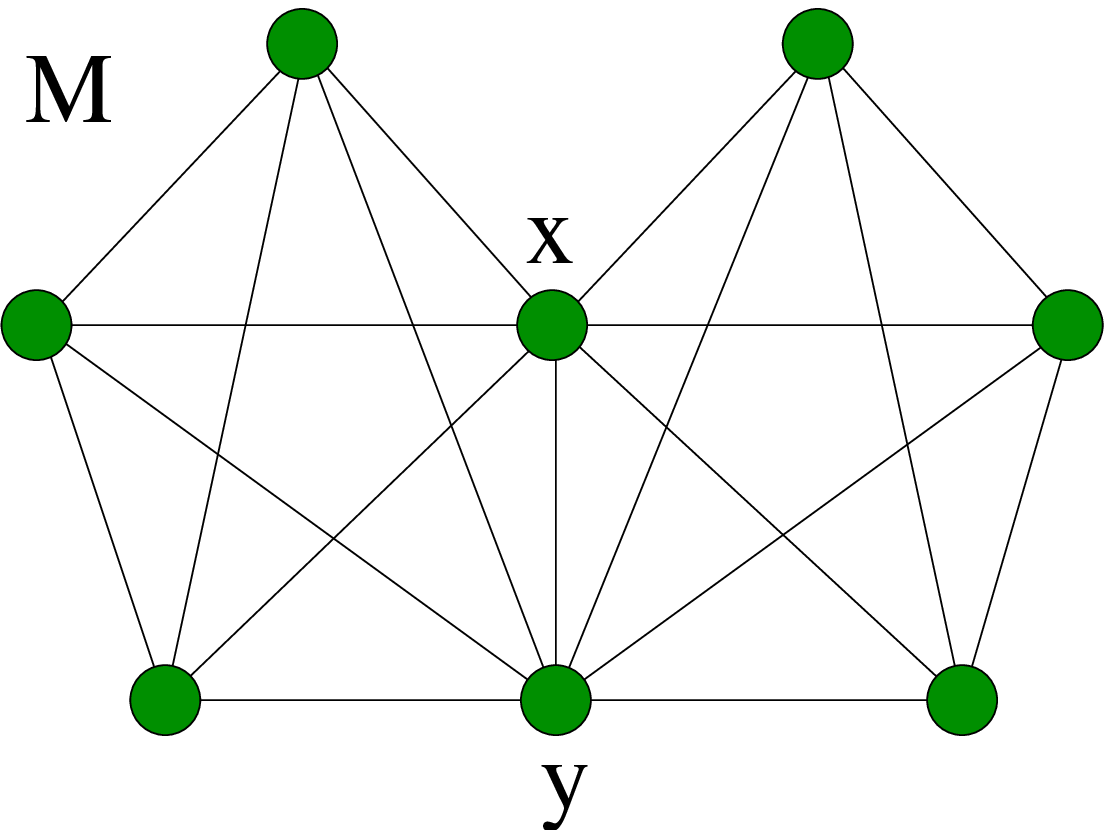}}
	\caption{$M$-graph.}
\end{figure}

The graph $K_5$ has six different embeddings in the torus as shown in Figure~3 (the torus is represented as a rectangle with opposite sides identified). To determine if a graph $G$ which has no $TK_{3,3}$'s embeds in the torus, it suffices to examine the ten
side components of $TK_5$ in $G$. By analyzing the six embeddings of $K_5$ on the torus, Gagarin and Kocay \cite{GK} have proven structural results
which imply the following theorem.

\begin{figure}[h]
	\centerline {\includegraphics[width=5.3in]{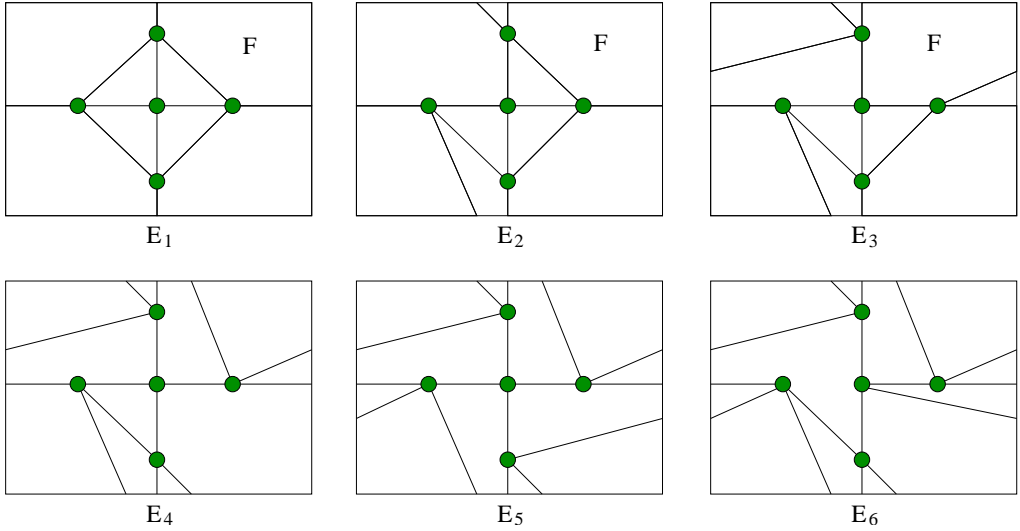}}
	\caption{Embeddings of $K_5$ in the torus.}
\end{figure}

\begin{thm}[Gagarin and Kocay \cite{GK}]
A graph $G$ with a $TK_5$ but no $TK_{3,3}$'s is toroidal if
and only if 

(i) all of the augmented side components of the $TK_5$ are planar, or

(ii) nine of the augmented side components of the $TK_5$
are planar and the remaining side component $S$ is special, or

(iii) $G$ contains a $TM$ where $M$ is as pictured in Figure~2 and
all of the augmented side components of the $TM$ are planar graphs.
\end{thm}
The special side component $S$ of Theorem 2(ii) must be embedded in a cylinder provided by the face $F$ of embeddings $E_1$ and $E_2$, or, in some particular cases, by the face $F$ of  embedding $E_3$, shown in Figure~3. Details can be found in \cite{GK}.

\section{Kuratowski theorem for the toroidal graphs with no $K_{3,3}$'s}
\noindent In this section we prove the characterization of the toroidal graphs with no $K_{3,3}$'s in terms of forbidden minors. Also we show the corresponding list of eleven topological obstructions.
\begin{thm}
A graph $G$ with no $K_{3,3}$'s is toroidal if and only if $G$ does not contain any of $G_1$, $G_2$, $G_3$ or $G_4$ of Figure~4 as a minor.
\end{thm}

\begin{figure}[h]
	\centerline {\includegraphics[width=4.2in]{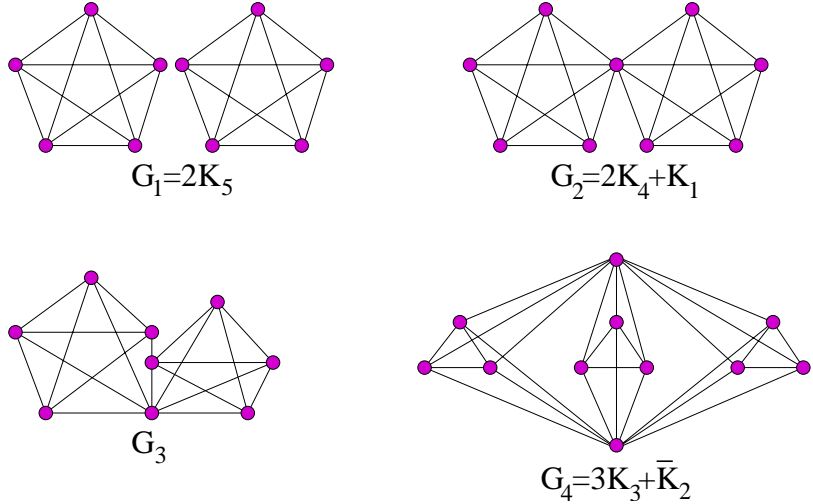}}
	\caption{\it The minor-minimal non-toroidal graphs $G_1, G_2, G_3$ and $G_4$ containing no $K_{3,3}$'s.}
\end{figure}
\begin{pf*}{Proof} \textsl{(Necessity).} The program from \cite{Wendy1} was used to check that none of the graphs $G_1, G_2, G_3$, or $G_4$ of Figure~4 is toroidal and also that deleting or contracting any edge of $G_1, G_2, G_3$, or $G_4$ results in a toroidal graph, i.e.  $G_1, G_2, G_3$ and $G_4$ are non-toroidal minor order minimal. This can also be verified by hand for $G_1, G_2$ and $G_3$ without too much difficulty (see Lemma 4.4 in \cite{GK}).

Graph $G_4$ of Figure~4 can be considered as obtained from the $M$-graph of Figure~2 by substituting side component $K_5-e$ for the central edge: the end-points $u$ and $v$ of the deleted edge $e=uv$ are identified respectively with the end-points $x$ and $y$ of the central edge $xy$ of the $M$-graph. Clearly, the corresponding augmented side component is $K_5$ which is non-planar, and, by Theorem\,2(iii), the whole graph $G_4$ is not toroidal. Deletion or contraction of any edge $e$ of $G_4$ results in an $M$-graph subdivision with all augmented side components planar. Therefore, by Theorem\,2(iii), $G_4$ is  a minor order minimal non-toroidal graph.

Since none of $G_1, G_2, G_3$, or $G_4$ is toroidal, it is obvious that any graph which contains one of them as a minor is not toroidal. Therefore we need to prove that these are the only minor order minimal obstructions.

\textsl{(Sufficiency).} If $v$ is a cut-vertex of $G$, then Miller \cite{Miller} has proven that the genus of $G$ equals the sum of the genera of the bridges of $G$ with respect to $v$. This implies that a graph $G$ is toroidal if and only if it contains at most one $2$-connected non-planar toroidal component and all the other $2$-connected components of $G$ are planar. Suppose $G$ does not contain any of $G_1, G_2, G_3, G_4$, or $K_{3,3}$ as a minor. Since $G$ has no $K_{3,3}$-minors, $G$ does not contain any $K_{3,3}$-subdivision as a subgraph. It is necessary to prove that $G$ is toroidal.

Since $G$ does not contain minors $G_1, G_2$, or $K_{3,3}$, $G$ can contain at most one non-planar $2$-connected component. If all the $2$-connected components of $G$ are planar, $G$ is toroidal. Therefore suppose $G$ is $2$-connected and non-planar. Since $G$ does not contain a $K_{3,3}$-subdivision, it contains a $K_5$-subdivision $TK_5$ such that (by Theorem\,1) each bridge of the set of corner vertices $K$ of $TK_5$ contains exactly two corner vertices from $K$, and it is possible to decompose $G$ into the side components of a $TK_5$.

Since $G$ does not contain minors $G_1$, or $G_2$, or $K_{3,3}$, at most one side component, say $H$, of $TK_5$ in $G$ can be non-planar. Figure~5 shows why there can be no pair of non-planar augmented side components of $TK_5$ in $G$: otherwise the side components $H_1, H_2$ and paths $ae, eb, de, ec$ in $G$ of Figure~5(a) would contain a forbidden minor $G_2$, or the side components $H_1, H_2$ and paths $bd, da, ae, ec$ in $G$ of Figure~5(b) would contain a forbidden minor $G_2$. The cases of other pairs of side components of $TK_5$ in $G$ are equivalent to the two cases of Figure~5 by symmetry. Therefore there is at most one non-planar augmented side component $H'$ of $TK_5$ in $G$. If $H'$ is planar, by Theorem\,2(i), $G$ is toroidal. Suppose $H'$ is the unique non-planar augmented side component of $TK_5$. Reasoning similar to that shown in Figure~5 demonstrates that the non-planar augmented side component $H'$ must be obtained from the planar or non-planar side component $H$. If $H$ is planar, by Theorem\,2(ii), $G$ is toroidal.

\begin{figure}[h]
	\centerline {\includegraphics[width=3.5in]{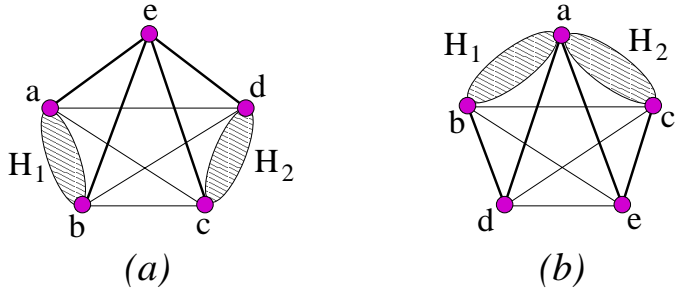}}
	\caption{\it Two side components $H_1$ and $H_2$ of $TK_5$ in $G$.}
\end{figure}
Suppose $H$ is the unique non-planar side component and $TK'_5$ is a $K_5$-subdivision in $H$. Because $G$ does not contain minors $G_1, G_2$ and $G_3$, the two subdivisions $TK_5$ and $TK'_5$ must share two common corners, say $x$ and $y$. Therefore $G$ contains a subdivision $TM$ of the $M$-graph (a detailed proof of this fact can be found in Theorem\,4.7 of \cite{GK}).

Now, because $G$ has no $K_{3,3}$-subdivision, we can decompose $G$ into the side components of the $M$-graph subdivision $TM$. Since $G$ does not contain minors $G_1, G_2$ and $K_{3,3}$, only the side component $H_{xy}$ of $TM$ in $G$ that corresponds to the central edge $xy$ of $M$-graph can be non-planar (see Figure~2): all the other side components of $TM$ in $G$ must be planar even when augmented.

Since $G$ does not contain minors $G_1, G_2, G_3$, or $K_{3,3}$, the augmented side component $H_{xy}+xy$ of $TM$ can contain just a $K_5$-subdivision $TK''_5$ that shares the same two common corners $x$ and $y$ with $TK_5$ and $TK'_5$ in $G$. However, because $G$ does not contain minor $G_4$, the side component $H_{xy}$ of $TM$ in $G$ must be planar even when augmented, and, by Theorem\,2(iii), the whole graph $G$ is toroidal.
\qed
\end{pf*}

The remaining seven topological obstructions that are not minor order obstructions are depicted in Figure~6. They were obtained by splitting the vertices of $G_1, G_2, G_3$, or $G_4$ in all possible ways by the first author, and were also determined as a part of a much larger project done by Chambers and Myrvold \cite{John}.\\
\begin{figure}[h!]
	\centerline {\includegraphics[width=4.1in]{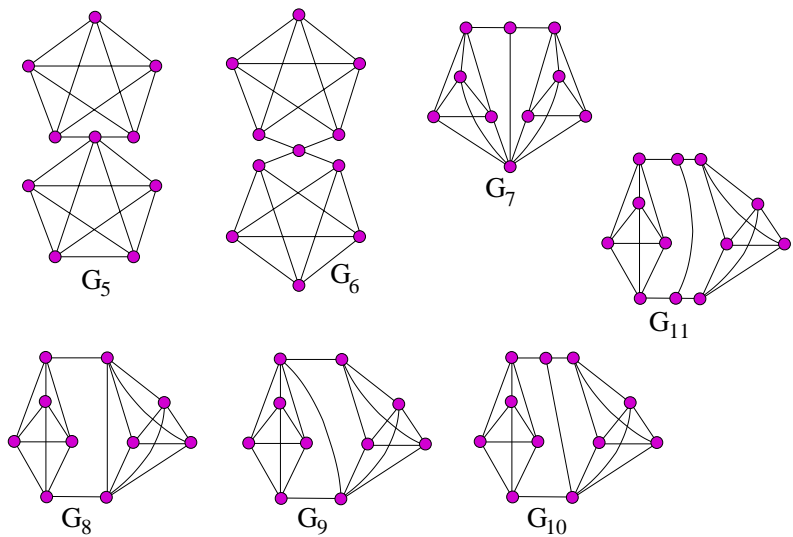}}
	\caption{\it Topological obstructions $G_5, G_6, \ldots, G_{11}$ for the torus.}
\end{figure}

\section{Final Remarks}
\noindent The results in this paper completely characterize the torus obstructions for all graphs with no $K_{3,3}$'s. Characterizing the complete set of torus obstructions seems to be a very difficult task. Such a characterization could help in finding a better understanding of the torus and, perhaps, the surfaces of higher genus. The obstructions for the toroidal graphs with no $K_{3,3}$'s suggest that a feasible approach to characterizing the complete set of torus obstructions would be to break the set into tractable subclasses.

Note that the torus can also be represented as a hexagon with its three pairs of opposite sides identified. This representation may help give a better understanding of the toroidal graphs as well. For example, the unique embedding of $K_7$ on the torus is very complicated in the rectangular representation of the torus. However there exists a simple symmetric drawing of $K_7$ on the torus represented as a hexagon (see Figure~7). The drawing of Figure~7 was found by the first author.\\
\begin{figure}[h!]
	\centerline {\includegraphics[width=1.7in]{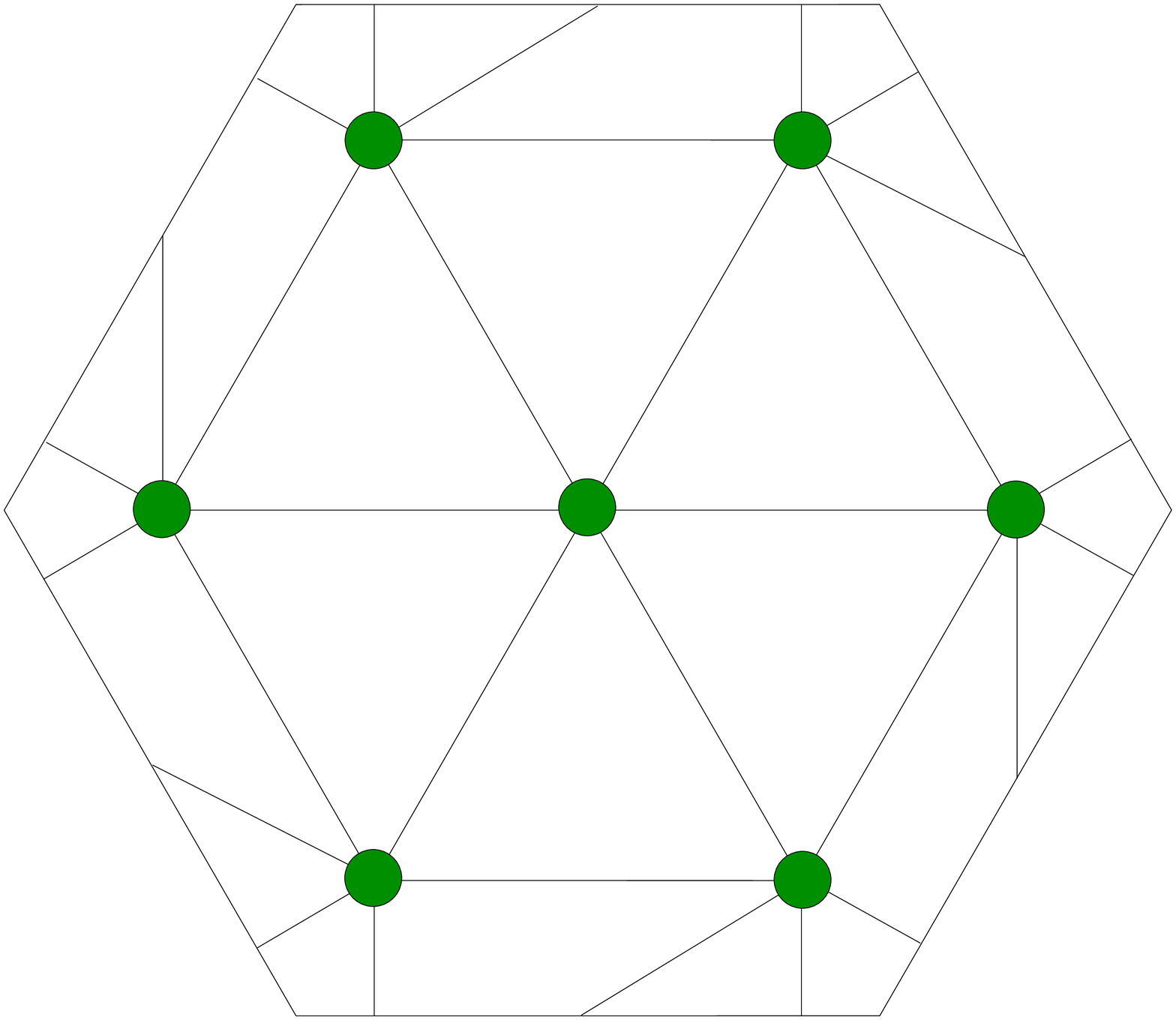}}
	\caption{\it A symmetric drawing of $K_7$ on the torus.}
\end{figure}

After the torus, the next most likely avenues of attack are determining complete obstruction sets for the Klein bottle, the spindle surface and the two-holed torus. Since a characterization in terms of forbidden minors and subdivisions can be too complicated, one may try to find a characterization of graphs on these surfaces in terms of forbidden substructures of a different kind.





\begin{thebibliography}{99}





\bibitem{DA81} ÊD. Archdeacon, A Kuratowski theorem for the projective plane, J. Graph Theory 5 (3) (1981) 243--246.

\bibitem{AB} D. Archdeacon and C. P. Bonnington, Obstructions for embedding cubic graphs on the spindle surface, J. Combin. Theory Ser. B, 91 (2) (2004) 229--252.

\bibitem{ArchHun} D. Archdeacon and J. P. Huneke, A Kuratowski Theorem for Nonorientable Surfaces, J. Combin. Theory Ser. B, 46 (1989) 173--231.

\bibitem{Asano} T. Asano, An approach to the subgraph homeomorphism problem, Theoret. Comput. Sci. 38 (1985) 249-267.

\bibitem{orient} R. Bodendiek and K. Wagner, Solution to {K}{\"{o}}nig's graph embedding problem, Math. Nachr. 140 (1989) 251-272.

\bibitem{Bondy} J.A. Bondy and U.S.R. Murty, Graph Theory
with Applications, American Elsevier Publishing, New York, 1976.

\bibitem{John} J. Chambers, Hunting for torus obstructions,
M.Sc. Thesis, Department of Computer Science, University of Victoria, 2002.

\bibitem{Diestel} R. Diestel, Graph Theory, 2nd edition,
Springer, New York, 2000.

\bibitem{Fellows} M. Fellows and P. Kaschube, Searching for
$K_{3,3}$ in linear time, Linear and Multilinear Algebra 29 (1991) 279-290.

\bibitem{GK} A. Gagarin and W. Kocay, Embedding 
graphs containing $K_5$-sub\-divisions, Ars Combin. 64 (2002) 33-49.

\bibitem{GLL2} A. Gagarin, G. Labelle, and P. Leroux, Characterization and enumeration of toroidal $K_{3,3}$-subdivision-free graphs, submitted, 2004. Preprint available at: {\texttt{http://arXiv.org/math.CO/0411356}}.

\bibitem{HG79} H. H. Glover, J. P. Huneke, and C. S. Wang, 103 graphs that are irreducible for the
		  projective plane, J. Combin. Theory Ser. B, 27 (3) (1979) 332-370.

\bibitem{Kelmans} A.K. Kelmans, Graph expansion and reduction, in: Algebraic methods in graph theory, Vol. I (Szeged, 1978), Colloq. Math. Soc. J\'anos Bolyai, 25, North-Holland, Amsterdam-New York, 1981, pp. 317--343.

\bibitem{Kuratowski} K. Kuratowski, Sur le probl\`eme des courbes
gauches en topologie, Fund. Math. 15 (1930), 271-283.

\bibitem{Miller} G. L. Miller, An Additivity Theorem for the Genus of a Graph, J. Combin. Theory Ser. B, 43 (1987) 25-47.

\bibitem{Wendy1} E. Neufeld and W. Myrvold, Practical Toroidality Testing, in:
Eighth Annual ACM-SIAM Symposium on Discrete Algorithms (SODA),
1997, pp. 574-580.

\bibitem{RS90} N. Robertson and P. D. Seymour, Graph minors. VIII. A Kuratowski theorem for general surfaces, J. Combin. Theory Ser. B, 48 (1990) 255-288.

\bibitem{Vollmer} H. Vollmerhaus, \"{U}ber die Einbettung von Graphen in zweidimensionale orientierbare Mannigfaltigkeiten kleinsten Geschlechts, in: Beitrage
z\"{u}r Graphentheorie, H. Sachs, H. Voss and H. Walther eds.,
Leipzig, 1968, pp. 163--168 (in German).

\bibitem{Wagner} K. Wagner, \"Uber eine Erweiterung eines Satzes von Kuratowski, Deutsche Math. 2 (1937) 280--285 (in German).

\bibitem{Wagner2} K. Wagner, \"{U}ber einer Eigenschaft der ebenen Komplexe, Math. Ann. 114 (1937) 570--590 (in German).

\end{thebibliography}
\end{document}